\documentclass{article}
\usepackage[margin=1.2in]{geometry}
\usepackage[english]{babel}
\usepackage{amsmath}
\usepackage{amsfonts}
\usepackage{amssymb}
\usepackage{amsthm}
\usepackage{tikz-cd}
\usetikzlibrary{positioning}
\usetikzlibrary{shapes.geometric}

\usepackage[colorlinks=true,breaklinks=true,bookmarks=true,urlcolor=blue,
     citecolor=blue,linkcolor=blue,bookmarksopen=false,draft=false]{hyperref}

\newcommand{\Z}		{\ensuremath{\mathbb{Z}}}

\newcommand{\R}		{\ensuremath{\mathbb{R}}}

\DeclareMathOperator{\Aut}{Aut}

\DeclareMathOperator{\im}{im}

\DeclareMathOperator{\ES}{ES}
\DeclareMathOperator{\coker}{coker}

\theoremstyle{plain}
\newtheorem{theorem}{Theorem}[section]

\newtheorem{lemma}[theorem]{Lemma}
\newtheorem{proposition}[theorem]{Proposition}

\theoremstyle{definition}
\newtheorem{definition}[theorem]{Definition}

\newtheorem{example}[theorem]{Example}

\makeatletter
\newcommand*{\doublerightarrow}[2]{\mathrel{
  \settowidth{\@tempdima}{$\scriptstyle#1$}
  \settowidth{\@tempdimb}{$\scriptstyle#2$}
  \ifdim\@tempdimb>\@tempdima \@tempdima=\@tempdimb\fi
  \mathop{\vcenter{
    \offinterlineskip\ialign{\hbox to\dimexpr\@tempdima+1em{##}\cr
    \rightarrowfill\cr\noalign{\kern.5ex}
    \rightarrowfill\cr}}}\limits^{\!#1}_{\!#2}}}
\makeatother

\title{Non-Abelian Flows in Networks}
\date{\today}
\author{D.M.H. van Gent\\ Leiden University}

\begin{document}
\maketitle
\begin{abstract}
In this work we consider a generalization of graph flows. A graph flow is, in its simplest formulation, a labeling of the directed edges with real numbers subject to various constraints. A common constraint is conservation in a vertex, meaning that the sum of the labels on the incoming edges of this vertex equals the sum of those on the outgoing edges. One easy fact is that if a flow is conserving in all but one vertex, then it is also conserving in the remaining one.
In our generalization we do not label the edges with real numbers, but with elements from an arbitrary group, where this fact becomes false in general.
As we will show, graphs with the property that conservation of a flow in all but one vertex implies conservation in all vertices are precisely the planar graphs.
\end{abstract}

\section{Introduction}

A {\em graph} (or {\em network}) is a pair \((V,E)\) where \(V\) is a finite set of {\em vertices} and the set of {\em edges} \(E\) is a subset of \(\binom{V}{2}\), the set of all size 2 subsets of \(V\).
In this article we consider groups which are not required to be abelian and therefore write our group operation multiplicatively.
With \(\Gamma\) a group and \(G=(V,E)\) a graph, we call a map \(f:V^2\to\Gamma\) a {\em\(\Gamma\)-flow in \(G\)} if for all \(u,v\in V\) we have \(f(u,v)=f(v,u)^{-1}\), and \(f(u,v)=1\) if \(\{u,v\}\not\in E\). 
This definition agrees with the classical definition of a network flow when \(\Gamma=\R\).

Non-abelian graph flows were first considered by M.J. DeVos in his PhD thesis \cite{DeVos} and later by A.J. Goodall et al. \cite{Goodall} and B. Litjens \cite{Litjens}.
They consider graphs embedded on surfaces and ask whether flows exists which are nowhere trivial, i.e. \(f(u,v)\neq 1\) if and only if \(\{u,v\}\in E\).
Although our main result involves planar embeddings of graphs, we instead ask to which extent Kirchhoff's law of conservation holds.

Let \(G=(V,E)\) be a graph, \(\Gamma\) a group and \(f\) a \(\Gamma\)-flow in \(G\).
We call \(f\) {\em tractable} if for each \(v\in V\) the subgroup \(\langle f(u,v) \,|\, u\in V \rangle\) of \(\Gamma\) is abelian.
For tractable \(f\) we define the {\em excess} \(e_f:V\to\Gamma\) to be the map given by \(v\mapsto \prod_{u\in V} f(u,v)\) and we say \(f\) is {\em conserving} in \(v\) if \(e_f(v)=1\).
In the classical case, we have the following lemma.

\begin{lemma}\label{lem:abelian_easy}
Let \(\Gamma\) be an {\em abelian} group, let \(f\) be a \(\Gamma\)-flow in a graph \(G=(V,E)\) and let \(w\in V\). If \(f\) is conserving in all vertices of \(V\setminus\{w\}\), then \(f\) is conserving in \(w\).
\end{lemma}
\proof
We have \(e_f(w) = \displaystyle\prod_{v\in V} e_f(v) = \displaystyle\prod_{(u,v)\in V^2} f(u,v) = \displaystyle\prod_{\{u,v\}\in E} f(u,v) f(v,u) = 1\). 
\endproof

We will show that Lemma~\ref{lem:abelian_easy} can fail for non-abelian \(\Gamma\). 
We say a flow \(f\) {\em leaks} if it is tractable and conserving in all but precisely one vertex and we call a graph \(G\) {\em leak-proof} if there exist no flows in \(G\) that leak for any group \(\Gamma\). 
Our main result, proven in Section~\ref{sec:planar_graphs}, is as follows.

\begin{theorem}\label{thm:main}
A graph is leak-proof if and only if it is planar.
\end{theorem}

We say a flow \(f\) of \(G\) has a {\em binary leak} at distinct vertices \(u,v\in V\) if it is tractable and conserving in all vertices of \(V\setminus\{u,v\}\) while \(e(u)e(v)\neq1\). Here \(u\) and \(v\) can be thought of as a source and sink of the flow. We call \(G\) {\em binary leak-proof} if no binary leaking flows exist for \(G\).
Analogously to Lemma~\ref{lem:abelian_easy} one can show that a flow cannot have a binary leak when the group is abelian. We also prove the following analogue to Theorem~\ref{thm:main} in Section~\ref{sec:extra_planar}.

\begin{definition}
We call a graph \(G=(V,E)\) {\em extra-planar} if for all pairs of distinct \(u,v\in V\) the graph \((V,E\cup\{u,v\})\) is planar.
\end{definition}

\begin{theorem}\label{thm:binleak}
A graph is binary leak-proof if and only if it is extra-planar.
\end{theorem}

Instead of studying leak-proof graphs, one could also study leak-proof groups, where we call a group \(\Gamma\) {\em leak-proof} if for all graphs \(G=(V,E)\) no tractable flows \(f:V^2\to\Gamma\) of \(G\) leak.
Theorem~\ref{thm:main} shows that the decision problem `Is this graph leak-proof?' can be decided in time \(O(|V|)\), as Hopcroft and Tarjan gave an algorithm to test graph planarity in \cite{PlanarityTesting} of this complexity. For leak-proof groups, we prove the following in Section~\ref{sec:group}.

\begin{theorem}\label{thm:lpgroup}
The decision problem `Is this finite group leak-proof?' is decidable.
\end{theorem}

The present work, in particular Theorem~\ref{thm:binleak}, was inspired by a problem the author encountered in his Master's thesis \cite{Thesis} on graded rings.
Here a flow with a binary leak gives rise to an example (Example~2.17 of \cite{Thesis}) of an efficient ring grading with a non-abelian group that cannot be replaced by an abelian group.

\section{Definitions and properties of (non-)planar graphs}

We briefly go through some basic definitions.
Let \(G=(V,E)\) be a graph.
We call a graph \((W,F)\) a {\em subgraph} of \(G\) if \(W\subseteq V\) and \(F\subseteq E\).
For \(W\subseteq V\) we call \((W,\{\{u,v\}\in E\,|\,u,v\in W\})\) the subgraph of \(G\) {\em induced by} \(W\).
With \(H=(W,F)\) a graph we call a map \(f:V\to W\) a {\em morphism} from \(G\) to \(H\) if \(f[E]\subseteq F\).
We call this \(f\) an {\em embedding} if it is injective and an {\em isomorphism} if \(f\) and its induced map \(E\to F\) are bijections.
A {\em path from \(u\in V\) to \(v\in V\)} in \(G\) is a finite sequence of vertices \((x_0,\dotsc,x_n)\) for some \(n\in\Z_{\geq 0}\) such that \(x_0=u\), \(x_n=v\) and \(\{x_i,x_{i+1}\}\in E\) for all \(0\leq i < n\). We call this path {\em non-trivial} if \(n>0\) and {\em closed} if \(x_0=x_n\).
We say \(u\in V\) is {\em connected} to \(v\in V\) in \(G\) if there exists a path from \(u\) to \(v\) in \(G\).
The `is connected to' relation is an equivalence relation on \(V\) and we call its equivalence classes the {\em connected components} of \(G\).
For \(u\in V\) we call \(v\in V\) a {\em neighbour} of \(u\) if \(\{u,v\}\in E\) and we write \(N_G(u)\subseteq V\) for the set of neighbours of \(u\).
An edge \(\{u,v\}\in E\) is called a {\em bridge} if all paths in \(G\) from \(u\) to \(v\) contain the edge \(\{u,v\}\).
A {\em forest} is a graph in which every edge is a bridge. 

We now give some facts about (non-)planar graphs.

\begin{definition}
For \(A,B\in\R^2\) write \(\overline{AB}\) for the line \(\{tA+(1-t)B\,|\,t\in(0,1)\}\).
Let \(G=(V,E)\) be a graph.
A {\em planar embedding} of \(G\) is an injective map \(\epsilon:V\to\R^2\) such that for all \(\{a,b\},\{c,d\}\in E\) we have \(\overline{\epsilon(a)\epsilon(b)} \cap \overline{\epsilon(c)\epsilon(d)}=\emptyset\) when \(\{a,b\}\neq\{c,d\}\), and \(\overline{\epsilon(a)\epsilon(b)}\cap \epsilon[V]=\emptyset\).
We call \(G\) {\em planar} if it has a planar embedding.
%Let \(\epsilon\) be a planar embedding of \(G\).
%A {\em face} of \((G,\epsilon)\) is a connected component of \(\R^2\setminus (\epsilon[V]\cup\bigcup_{\{u,v\}\in E} \overline{\epsilon(u)\epsilon(v)})\).
%We say \(\{a,b\}\in E\) {\em borders} a face \(F\) of \((G,\epsilon)\) if \(\overline{\epsilon(a)\epsilon(b)}\) is contained in the topological closure of \(F\).
\end{definition}

The above definition of a planar embedding has been simplified for our purposes, which is justified by F\'ary's Theorem \cite{Fary}.

%Note that each edge borders at most two faces of any planar embedding, namely the faces `on either side' of the edge.
%These faces can however be the same.

\begin{definition}\label{def:orientation}
Let \(G=(V,E)\) be a graph with a planar embedding \(\epsilon\).
The {\em orientation} of \((G,\epsilon)\) at \(v\in V\) is the clockwise permutation \(\rho_\epsilon(v)\) of \(N_G(v)\).
A {\em boundary walk} of \((G,\epsilon)\) is a non-trivial closed path \((x_0,x_2,\dotsc,x_n)\) in \(G\) such that for all \(i,j\in\Z/n\Z\) we have \(x_{i+2}=\rho_\epsilon(x_{i+1})(x_i)\) and if \((x_i,x_{i+1})=(x_j,x_{j+1})\), then \(i=j\).
\end{definition}

\begin{lemma}\label{lem:boundary_walk}
Let \(\epsilon\) be a planar embedding of a graph \(G=(V,E)\) and let \(p=(u_1,u_2,\dotsc,u_n)\) be a boundary walk.
If \((u_i,u_{i+1})=(u_{j+1},u_{j})\) for some \(i,j\in\Z/n\Z\), then \(\{u_i,u_{j}\}\) is a bridge.
\end{lemma}
\proof
To show that \(e=\{u_i,u_{j}\}\) is a bridge, it suffices to show that \(u_i\) and \(u_j\) are disconnected in the graph \(G'=(V,E')\) with \(E'=E\setminus\{e\}\).
Note that \(a,b\in V\) are connected in \(G'\) if and only if \(\epsilon(a)\) and \(\epsilon(b)\) are connected in the topological space \(X=\epsilon[V]\cup\bigcup_{\{x,y\}\in E'}\overline{\epsilon(x)\epsilon(y)}\). 
Hence it suffices by the Jordan curve theorem to show that there exists a loop \(C\) in \(\R^2\setminus X\) separating \(u_i\) and \(u_j\), as any path from \(u_i\) to \(u_j\) must intersect this loop.

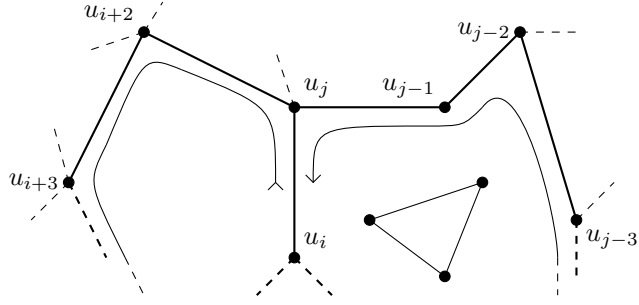
\begin{figure}[h]
\centering
\begin{tikzpicture}
\draw[fill=black] (0,0) circle[radius=2pt] node[above right] {$u_i$};
\draw[fill=black] (0,2) circle[radius=2pt] node[above right] {$u_{j}$};
\draw[fill=black] (-2,3) circle[radius=2pt] node[above left] {$u_{i+2}$};
\draw[fill=black] (-3,1) circle[radius=2pt] node[left] {$u_{i+3}$};
\draw[fill=black] (2,2) circle[radius=2pt] node[above left] {$u_{j-1}$};
\draw[fill=black] (3,3) circle[radius=2pt] node[left] {$u_{j-2}$};
\draw[fill=black] (3.75,.5) circle[radius=2pt] node[below right] {$u_{j-3}$};
\draw[fill=black] (2.5,1) circle[radius=2pt] node[below right] {};
\draw[fill=black] (1,.5) circle[radius=2pt] node[below right] {};
\draw[fill=black] (2,-.25) circle[radius=2pt] node[below right] {};
\draw (2.5,1) -- (1,.5) -- (2,-.25) -- cycle;
\draw[thick] (0,0) -- (0,2) -- (-2,3) -- (-3,1);
\draw[thick] (0,2) -- (2,2) -- (3,3) -- (3.75,.5);
\draw[thick,dashed] (0,0) -- (-.5,-.5);
\draw[thick,dashed] (0,0) -- (.5,-.5);
\draw[thick,dashed] (3.75,.5) -- (3.75,-.25);
\draw[thick,dashed] (-3,1) -- (-2.5,0);
\draw[dashed] (0,2) -- (-.25,2.75);
\draw[dashed] (-2,3) -- (-2.75,2.75);
\draw[dashed] (-2,3) -- (-1.75,3.4);
\draw[dashed] (-3,1) -- (-3.5,.5);
\draw[dashed] (-3,1) -- (-3.2,1.75);
\draw[dashed] (3,3) -- (3.75,3);
\draw[dashed] (3.75,.5) -- (4.25,1);
\draw (-.35,.9) -- (-.25,1) -- (-.15,.9); 
\draw (-.25,1) .. controls (-.25,2) .. (-1.5,2.5);
\draw (-1.5,2.5) .. controls (-2,2.7) .. (-2.5,1.5);
\draw (-2.5,1.5) .. controls (-2.75,1) .. (-2.25,0);
\draw[dashed] (-2.25,0) -- (-2,-.5);
\draw (.35,1.1) -- (.25,1) -- (.15,1.1); 
\draw (.25,1) .. controls (.25,1.75) .. (1.75,1.75);
\draw (1.75,1.75) .. controls (2.25,1.75) .. (2.5,2);
\draw (2.5,2) .. controls (3,2.65) and (3.5,.5) .. (3.5,0);
\draw[dashed] (3.5,0) -- (3.5,-.5);
%\draw (3.5,.5) .. controls (3.75,-.25) .. (4,.5);
%controls (2,2) .. (2.5,2.5);
\end{tikzpicture}
\caption{Boundary walk}
\end{figure}
We informally construct this loop as follows (see Figure 1).
Place yourself at the midway point between \(u_i\) and \(u_j\).
Walk along the path \(p\) in \(G\) in the direction of \(u_j\) and while doing so draw a continuous curve \(C\) on your left hand side, being careful not to let \(C\) intersect itself or the graph.
That this is possible follows from the definition of a boundary walk.
Stop once you have reached your starting point for the first time again, and note that this time you are facing \(u_i\) by the assumption that \((u_i,u_{i+1})=(u_{j+1},u_j)\).
Thus on your right hand side is the start of your curve \(C\), and connect the endpoints, crossing \(\overline{\epsilon(u_i)\epsilon(u_j)}\) once.
Then \(C\) satisfies the requirements, so \(e\) is a bridge.
\endproof

\begin{definition}\label{def:contraction}
Let \(G=(V,E)\) be a graph.
We call a subgraph \(H=(W,F)\) of \(G\) a {\em spanning forest} if it is a forest and \(W=V\).
For a spanning forest \(H=(W,F)\) of \(G\) we define \(G_H=(C,D)\) to be the {\em contraction of \(H\) in \(G\)}, where \(C\) is the set of connected components of \(H\) and \(D=\{\{X,Y\}\in\binom{C}{2}\,|\,(\exists\, u\in X,\,v\in Y) \ \{u,v\}\in E\}\).
A graph \(M\) is a {\em minor} of \(G\) if it can be embedded in some contraction of \(G\).
\end{definition}

Note that the `is a minor of' relation is a partial order (up to graph isomorphism). In particular, if \(I\) is a minor of \(H\) and \(H\) is a minor of \(G\), then \(I\) is a minor of \(G\).
Write \(K_5\) for the complete graph on 5 vertices and \(K_{3,3}\) for the complete bipartite graph on 3 and 3 vertices.

\begin{theorem}[Kuratowski, Theorem 4.4.6 in \cite{gtbook}]\label{Kuratowski}
A graph \(G\) is planar if and only if \(G\) does not have \(K_5\) or \(K_{3,3}\) as a minor.
\end{theorem}

\section{Non-planar graphs}

First we show that all non-planar graphs can leak.

\begin{lemma}\label{lem:subgraph_leak_proof}
A graph is leak-proof if and only if all its subgraphs are leak-proof.
\end{lemma}
\proof
Since each graph is its own subgraph, the implication (\(\Leftarrow\)) is trivial.
Let \(G=(V,E)\) be a graph with a subgraph \(H=(W,F)\) and assume that there exists some group \(\Gamma\) with a leaking \(\Gamma\)-flow \(g:W^2\to\Gamma\) of \(H\). 
Then we consider \(f:V^2\to\Gamma\) by taking \(f(u,v)=g(u,v)\) when \(\{u,v\}\in F\) and \(f(u,v)=1\) otherwise. Then \(f\) is a leaking flow for \(G\), proving (\(\Rightarrow\)).
\endproof

\begin{proposition}\label{prop:minor_leak}
A graph is leak-proof if and only if all its minors are leak-proof.
\end{proposition}
\proof
Let \(G=(V,E)\) be a graph. By Lemma~\ref{lem:subgraph_leak_proof} it suffices to show that if a contraction of a spanning tree \(H\) in \(G\) admits a leaking flow, then so does \(G\).
By induction we may even assume \(H\) has only a single edge \(e=\{a,b\}\). Then \(G_H\cong(W,F)\) with \(W=(V\setminus e)\cup\{e\}\) under the natural isomorphism \(e\mapsto e\) and \(w\mapsto\{w\}\) for \(w\in V\setminus e\). Assume \((W,F)\) admits a flow \(f:W^2\to\Gamma\) leaking at \(w\in W\) for some group \(\Gamma\).
Let \(X=N_G(a)\setminus e\) and \(Y=N_G(b)\setminus (e\cup X)\).
We define a flow \(g:V^2\to\Gamma\) such that for \(u,v\in W\) it is given by 
\begin{align*}
g(u,v)&=f(u,v) &u,v\not\in e, \\
g(a,u)^{-1}=g(u,a)&=f(u,e) & u\in X, \\
g(v,b)^{-1}=g(b,v)&=f(e,v) & v\in Y, \\
g(b,a)^{-1}=g(a,b)&=\prod_{u\in X \setminus\{b\}} f(u,a),
\end{align*}
and \(g(u,v)=1\) otherwise. Note that \(g\) agrees with \(f\) outside of \(e\) and that the flows going to \(e\) have been divided among \(a\) and \(b\).
Thus \(g\) is tractable and \(e_g(u)=e_f(u)\) for \(u\not\in e\). By definition of \(g(a,b)\) we have that \(e_g(a)=1\) and \(e_g(b)=e_f(e)\).
Hence \(g\) is a leaking flow for \(G\).
\endproof

To show that non-planar graphs are not leak-proof, it now suffices by Theorem~\ref{Kuratowski} to show that \(K_5\) and \(K_{3,3}\) admit a leaking flow.

\begin{definition}
Let \(C_2\) be the cyclic group with two elements.
Let \(n\in\Z_{>0}\) and consider the groups \(N=C_2^{n+1}=\langle z, x_1,\dotsc,x_n\rangle\) and \(G=C_2^n=\langle x_{n+1},\dotsc,x_{2n}\rangle\).
Define an action \(\varphi:G\to\Aut(N)\) defined on the generators as
\[ x_{n+i} \mapsto \big( x_j \mapsto x_j z^{\delta_{ij}},\quad z \mapsto z \big) \quad \text{for all } 1\leq i,j \leq n, \]
where \(\delta_{ij}=1\) if \(i=j\) and \(\delta_{ij}=0\) otherwise.
Then define the group \(\ES_n=N \rtimes_\varphi G\).
\end{definition}

Although we will not use the fact, the \(\ES_n\) are all {\em extraspecial \(2\)-groups}.

\begin{example}\label{example:k33}
Consider the utility graph \(K_{3,3}=(V,E)\) with \(V=\{1,2,3,4,5,6\}\) and \(E=\{\{u,v\}\,|\,u\in\{1,2,3\},\,v\in\{4,5,6\}\}\).
We define a flow \(f:V^2\to \ES_2\) which we specify by an \(\ES_2\)-valued (symmetric) matrix where the omitted entries are trivial:
\begin{align*}
f = \left(\begin{array}{ccc|ccc}
 & & & x_1 & x_2 & x_1 x_2 \\
 & & & x_4 & x_3 & x_4 x_3 \\
 & & & x_1 x_4 & x_2 x_3 & x_1 x_4 x_2 x_3 \\ \hline
x_1 & x_4 & x_1 x_4 & & & \\
x_2 & x_3 & x_2 x_3 & & & \\
x_1 x_2 & x_4 x_3 & x_1 x_4 x_2 x_3 & & &
\end{array}\right).
\end{align*}
For the first 5 columns it is easy to see that multiplying the first two non-trivial entries yields the third.
Thus for the first five vertices \(v\) we have \(\langle f(u,v) \,|\, u\in V \rangle\cong C_2^2\), which is abelian, and \(e_f(v)=1\).
For \(v=6\) we observe that \((x_1x_2)(x_4x_3)(x_1x_4x_2x_3)=z\) and thus \(\langle f(u,6)\,|\, u\in V\rangle=\langle x_1x_2,x_4x_3,z\rangle\cong C_2^3\) is abelian, and \(e_f(6)=z\neq 1\).
Hence \(f\) is a tractable flow that leaks at \(6\) and \(K_{3,3}\) is not leak-proof.
\end{example}

\begin{example}\label{example:k5}
Consider the complete graph \(K_{5}=(V,E)\) with \(V=\{1,2,3,4,5\}\).
Now we consider \(f:V^2\to \ES_3\) given by 
\begin{align*}
f = \left(\begin{array}{ccccc}
 & x_1 & x_2 & x_3 & x_1x_2x_3 \\
x_1 & & x_6 & x_5 & x_1 x_6 x_5 \\
x_2 & x_6 & & x_4 & x_2x_6x_4 \\
x_3 & x_5 & x_4 & & x_3x_5x_4 \\
x_1x_2x_3 & x_1x_6x_5 & x_2x_6x_4 & x_3x_5x_4 &
\end{array}\right).
\end{align*}
For each of the first four columns one notes that its first three non-trivial elements commute pair-wise, while multiplying them yields the fourth.
Thus for the first four vertices \(v\) the group \(\langle f(u,v)\,|\, u\in V\rangle\cong C_2^3\) is abelian and \(e_f(v)=1\).
For the last column, note that each pair \((a,b)\) of entries is of the form \(a=x_ix_jx_k\) and \(b=x_ix_{j+3}x_{k+3}\) with \(i,j,k,j+3,k+3\in\Z/6\Z\) distinct.
Hence \(ab=x_i^2 (x_j x_k) (x_{j+3} x_{k+3}) = x_i^2 (x_{j+3}x_{k+3})(x_j x_k) = ba \), so each pair commutes.
Finally, one computes \(e_f(5)=(x_1x_2x_3)(x_1x_6x_5)(x_2x_6x_4)(x_3x_5x_4)=z\neq 1\).
Thus \(f\) is a tractable leaking flow and thus \(K_5\) is not leak-proof.
\end{example}

Both examples were found by starting with the free group \(F\) with symbols \(V^2\) and dividing out the relations \(N\trianglelefteq F\) required to make the obvious map \(f:V^2\to F/N\) a tractable flow that is conserving in \(\#V-1\) vertices. Adding the restriction that the generators have order \(2\) gives us the groups \(\ES_2\) and \(\ES_3\). 

It now follows that all non-planar graphs leak, so we are half-way done proving Theorem~\ref{thm:main}.

\section{Planar graphs} \label{sec:planar_graphs}

Now we will prove that all planar graphs are leak-proof by induction.
For this we require a definition of the excess for non-tractable flows.

\begin{definition}
Let \(G=(V,E)\) be a graph with planar embedding \(\epsilon\) and let \(f:V^2\to\Gamma\) be a flow of \(G\).
Write \(C(\Gamma)\) for the set of conjugacy classes of \(\Gamma\) and \(\equiv\) for equality up to conjugation.
Then we define for \((G,\epsilon,f)\) the {\em round flow} \(r:V\to C(\Gamma)\) as \(r(v)\equiv 1\) if \(N_G(V)=\emptyset\), and otherwise
\[r(v)\equiv f(\rho_\epsilon(v)^0(u),v)\cdot f(\rho_\epsilon(v)^1(u),v)\dotsm f(\rho_\epsilon(v)^{n-1}(u),v),\]
where \(u\in N_G(v)\), \(n=\# N_G(v)\) and \(\rho_\epsilon\) is as in Definition~\ref{def:orientation}.
\end{definition}

Note that choosing a different \(u\in N_G(v)\) in the above definition results in a cyclic permutation of the factors, hence the products are conjugate in \(\Gamma\).
Thus the round flow is well-defined. Since \(1\in\Gamma\) is only conjugate to itself, we have that \(r(v)\equiv 1\) if and only if \(e(v)=1\) when the latter is defined.

\begin{theorem}\label{theorem:planar_cycle_stability}
Let \(G=(V,E)\) be a graph with planar embedding \(\epsilon\) and let \(f:V^2\to\Gamma\) be a flow of \(G\).
Let \(u\in V\) and assume \(r(v)\equiv1\) for all \(v\in V\setminus\{u\}\).
Then \(r(u)\equiv 1\).
\end{theorem}
\proof
Firstly, if \(G\) is the singleton graph, then \(r(u)\equiv 1\) is the empty product, so we are done.
We now apply induction and thus assume that the statement holds for all strict subgraphs \((W,F)\) of \(G\) with planar embedding \(\epsilon|_W\).

Secondly, we consider the case where \(G\) is not connected. 
Here we may apply the induction hypothesis to the induced subgraph of \(G\) with as vertex set the connected component of \(u\) to conclude that \(r(u)\equiv 1\). 

Thirdly, we consider the case where \(G\) is a forest. 
Then \(G\) has at least two vertices of degree 1, of which one, say \(v\), is not \(u\).
Let \(e=\{v,w\}\in E\) be the unique edge incident to \(v\), and note that \(1\equiv r_f(v)\equiv f(w,v)\) implies \(f(w,v)=1\). 
Hence \(f\) is a flow of the subgraph \(H\) of \(G\) obtained by removing \(e\).
Note that \(\epsilon\) is a planar embedding of \(H\) with the same round flow in each vertex, hence by the induction hypothesis we have \(r_f(u)\equiv1\).

Lastly we consider the case where \(G\) is connected and not a forest.
Then \(G\) has an edge \(\{v,w\}\in E\) that is not a bridge.
Then by Lemma~\ref{lem:boundary_walk} the boundary walk \(p=(x_0,\dotsc,x_n)\) of \((G,\epsilon)\) with \(x_0=v\) and \(x_1=w\) satisfies \((w,v)\neq(x_i,x_{i+1})\) for all \(i\in\Z/n\Z\).
Let \(b:V^2\to\{0,1\}\) be the map such that for all \(s,t\in V\) we have \(b(s,t)=1\) if and only if there exists some \(i\in\Z/n\Z\) such that \((s,t)=(x_i,x_{i+1})\).
Now consider \(\gamma=f(v,w)\) and \(g:V^2\to \Gamma\) given by
\[ (s,t) \mapsto \gamma^{b(t,s)} \cdot f(s,t) \cdot \gamma^{-b(s,t)}. \]
Firstly note that \(g\) is a flow of \(G\):
For all \(s,t\in V\) we have 
\[g(s,t)^{-1}=\gamma^{b(s,t)}\cdot f(s,t)^{-1}\cdot \gamma^{-b(t,s)}=g(t,s)\]
since \(f\) is a flow, and if \(\{s,t\}\not\in E\) we have \(g(s,t)=f(s,t)=1\) as \(b(s,t)=b(t,s)=0\).
Secondly, we have that \(g(v,w)=\gamma^{0} \cdot \gamma \cdot \gamma^{-1}=1\) by choice of \(\{v,w\}\), so \(g\) is even a flow of the subgraph \(H\) of \(G\) obtained by removing \(\{v,w\}\).
We now show that the round flows \(r_f\) and \(r_g\) of \(f\) respectively \(g\) in \((G,\epsilon)\) are conjugates in \(\Gamma\) at each vertex.
Then by the induction hypothesis applied to \(H\) it follows that \(r(u)\equiv 1\). 
Note that for all \(s,t\in V\) we have by definition of \(b\) that \(b(t,s)=b(s,\rho_\epsilon(s)(t))\).
Using this, we now simply verify for \(\{s,t\}\in E\), \(n=\# N_G(s)\) and \(\rho=\rho_\epsilon(s)\) that
\begin{align*}
r_{g}(s) 
&\equiv \prod_{k=0}^{n-1} g(\rho^k(t),s) \equiv \prod_{ k=0 }^{n-1} \gamma^{b(s,\rho^k(t))} \cdot f(\rho^k(t),s ) \cdot  \gamma^{-b(\rho^k(t),s)} \\
&\equiv \gamma^{b(s,t)} \left( \prod_{k=0}^{n-1} f(\rho^k(t),s) \gamma^{-b(\rho^k(t),s)} \gamma^{b(s,\rho^{k+1}(t))}  \right) \gamma^{-b(s,\rho^n(t))} \\
&\equiv \gamma^{b(s,t)} \left( \prod_{k=0}^{n-1} f(\rho^k(t),s) \right) \gamma^{-b(s,t)} \equiv \prod_{k=0}^{n-1} f(\rho_v^k(t),s ) \equiv r_f(s),
\end{align*}
as was to be shown.
We conclude that the statement holds for all planar graphs by induction.
\endproof

An earlier proof of Theorem~\ref{theorem:planar_cycle_stability} was due to H.W. Lenstra. In his version he does not remove edges in the inductive step but  contracts them in the sense of Definition~\ref{def:contraction}. This proof turned out to be more difficult to formalize. 

\proof{Proof of Theorem~\ref{thm:main}.} 
A non-planar graph has either \(K_5\) or \(K_{3,3}\) as minor by Theorem~\ref{Kuratowski}.
Both \(K_5\) and \(K_{3,3}\) are not leak-proof by Example~\ref{example:k5} respectively Example~\ref{example:k33}, so by Proposition~\ref{prop:minor_leak} neither are the non-planar graphs.
Let \(G=(V,E)\) be a planar graph with \(u\in V\) and let \(f\) be a tractable flow of \(G\) such that \(e(v)=1\) for all \(v\in V\setminus\{u\}\).
After choosing a planar embedding for \(G\) we have \(r(u)\equiv 1\) by Theorem~\ref{theorem:planar_cycle_stability} and thus \(e(u)=1\).
Hence \(f\) does not leak and \(G\) is leak-proof. 
\endproof

\section{Extra-planar graphs}\label{sec:extra_planar}

In this section we will prove Theorem~\ref{thm:binleak}, classifying the binary leak-proof graphs.
To do this we first prove a `Kuratowski's Theorem' for extra-planar graphs. 
Write \(K_5^{-}\) and \(K_{3,3}^{-}\) for the graphs obtained from \(K_5\) respectively \(K_{3,3}\) by removing a single edge, which by symmetry we do not have to specify.

\begin{theorem}\label{thm:extra_kuratowski}
A graph \(G\) is extra-planar if and only if \(G\) does not have \(K_5^-\) or \(K_{3,3}^-\) as a minor.
\end{theorem}
\proof
(\(\Rightarrow\)) This follows directly from Kuratowski's Theorem: If \(K_5^{-}\) or \(K_{3,3}^{-}\) is a minor of \(G\), then we may add a single edge to \(G\) such that \(K_5\) respectively \(K_{3,3}\) becomes a minor of this new graph, which is then non-planar.

\((\Leftarrow)\) We proceed by contraposition, so assume that \(G\) is not extra-planar. 
Let \(u,v\in V\) be such that \(G^+=(V,E\cup\{\{u,v\}\})\) is non-planar and let \(H^+=(V,F)\) be a spanning forest of \(G^+\) such that \(K_5\) or \(K_{3,3}\) embeds into \(G^+_H\). 
%If \(\{u,v\}\not\in H^+\) then \(H^+\) is a spanning forest of \(G\) and \(G_{H^+}\) equals \(G^+_{H^+}\) missing possible a single edge between the connect components of \(u\) and \(v\).
%In this case \(G_{H^+}\) clearly has \(K_5^-\) or \(K_{3,3}^-\) as minor.
Consider the spanning forest \(H=(V,F\setminus\{\{u,v\}\})\) of \(G\).
Then \(H\) has the same connected components as \(H^+\) with the exception that if \(H^+\) has a connected component containing both \(u\) and \(v\), it might have been split into two.  
Let \(T_u\) and \(T_v\) be the connected components of \(u\) respectively \(v\) in \(H\).

\begin{figure}[h]
\centering
\begin{minipage}{.5\textwidth}
  \centering
	\begin{tikzpicture}
	\node (U) at (0,0) {\(T_u\)};
	\node (V) [below=.5cm of U] {\(T_v\)};
	\node (A) [right=2.5cm of V] {\(R_1\)};
	\node (B) [below=1cm of A] {\(R_2\)};
	\node (C) [below=1cm of B] {\(R_3\)};
	\node (D) [below=1cm of V] {\(L_2\)};
	\node (E) [below=1cm of D] {\(L_3\)};
	\draw [semithick,dashed] (U) -- (A);
	\draw [semithick,dashed] (U) -- (B);
	\draw [semithick,dashed] (U) -- (C);
	\draw [semithick,dashed] (V) -- (A);
	\draw [semithick,dashed] (V) -- (B);
	\draw [semithick,dashed] (V) -- (C);
	\draw [semithick] (D) -- (A);
	\draw [semithick] (D) -- (B);
	\draw [semithick] (D) -- (C);
	\draw [semithick] (E) -- (A);
	\draw [semithick] (E) -- (B);
	\draw [semithick] (E) -- (C);
	\end{tikzpicture}
	\caption{Case \(K_{3,3}\)}
\end{minipage}%
\begin{minipage}{.5\textwidth}
  \centering
	\begin{tikzpicture}
	\def\ngon{6}
	\def\mylist{\(T_v\),\(T_u\),\(P_1\),\(P_2\),\(P_3\),\(P_4\)};
	\node[regular polygon,regular polygon sides=\ngon,minimum size=4.3cm] (p) {};
	\foreach[count=\x] \s in \mylist {\node (p\x) at (p.corner \x) {\s};}
	\foreach\x in {3,...,\numexpr\ngon-1\relax}{
	  \foreach\y in {\x,...,\ngon}{
	    \draw [semithick] (p\x) -- (p\y);
	  }
	}
	\foreach\x in {1,...,2}{
	  \foreach\y in {3,...,\numexpr\ngon\relax}{
	    \draw [semithick,dashed] (p\x) -- (p\y);
	  }
	}
	\end{tikzpicture}
  \caption{Case \(K_5\)}
\end{minipage}
\end{figure}

{\em Case \(K_{3,3}\):} First consider the case where \(K_{3,3}\) embeds into \(G^+_{H^+}\), meaning there is a subset \(C=\{L_1,L_2,L_3,R_1,R_2,R_3\}\) of size 6 of the set of connected components of \(H^+\) such that \(S^+=(C,\{\{L_i,R_j\} \,|\, i,j\in\{1,2,3\}\})\) is a subgraph of \(G^+_{H^+}\). 
If all elements of \(C\) are also connected components of \(H\), then \(G_H\) has the graph \(S^+\) minus possibly a single edge induced by \(\{u,v\}\) as subgraph, hence \(G\) has \(K_{3,3}^-\) as a minor.
Otherwise, for some \(X\in C\) we have \(X=T_u\sqcup T_v\) and without loss of generality \(X=L_1\).
Then the subgraph \(S\) of \(G_H\) induced by \(\{T_u,T_v,L_2,L_3,R_1,R_2,R_3\}\) is as in Figure 2, where the dashed lines indicate edges which are possibly present.
Merging \(T_u\) and \(T_v\) in \(S\) yields \(S^+\cong K_{3,3}\), hence for each \(i\in\{1,2,3\}\) the edge \(\{T_u,R_i\}\) or \(\{T_v,R_i\}\) is present.
Thus \(T_u\) or \(T_v\) has degree at least \(2\), which without loss of generality is \(T_v\).
It follows that \(K_{3,3}^-\) embeds into the subgraph of \(G_H\) induced by \(\{T_v,L_2,L_3,R_1,R_2,R_3\}\), so \(K_{3,3}^-\) is a minor of \(G\).

{\em Case \(K_5\):}
Now consider the case \(K_5\) embeds into \(G^+_{H^+}\), meaning there is a subset \(C=\{P_1,\dotsc,P_5\}\) of the set of connected components of \(H^+\) such that the subgraph of \(G^+_{H^+}\) induced by \(C\) is isomorphic to \(K_5\).
As before, the only interesting case is where \(P_5=T_u\sqcup T_v\).
Then the subgraph \(S\) of \(G_H\) induced by \(\{T_u,T_v,P_1,P_2,P_3,P_4\}\) is as in Figure 3.
Since merging \(T_u\) and \(T_v\) in \(S\) yields \(K_{5}\), for each \(i\in\{1,\dotsc,4\}\) the edge \(\{T_u,P_i\}\) or \(\{T_v,P_i\}\) is present.
If both \(T_u\) and \(T_v\) have degree \(2\), then without loss of generality \(S\) contains the edges \(\{T_u,P_3\}\), \(\{T_u,P_4\}\), \(\{T_v,P_1\}\) and \(\{T_v,P_2\}\).
Now note that \(S\) contains a \(K_{3,3}^-\) which partitions its vertices as \(\{\{T_u,P_1,P_2\},\{T_v,P_3,P_4\}\}\).
Hence \(G\) contains \(K_{3,3}^-\) as a minor.
Otherwise, without loss of generality \(T_u\) has degree at least \(3\) in \(S\) and the subgraph of \(G_H\) induced by \(\{T_u,P_1,\dotsc,P_4\}\) is either \(K_5\) or \(K_5^-\).
Hence \(G\) has \(K_5^-\) as a minor.

As \(G\) has \(K_{3,3}^-\) or \(K_5^-\) as a minor, the claim follows.
\endproof

We are now able to prove Theorem~\ref{thm:binleak}.

\proof{Proof of Theorem~\ref{thm:binleak}.} 
(\(\Leftarrow\))
Let \(G=(V,E)\) be an extra-planar graph and let \(f:V^2\to\Gamma\) be a tractable flow of \(G\) such that there are distinct \(u,v\in V\) with \(e_f(w)=1\) for all \(w\in V\setminus\{u,v\}\). Consider the graph \(H=(V,E\cup\{\{u,v\}\})\) and let \(\epsilon\) be a planar embedding of \(H\).
Now let \(g:V^2\to\Gamma\) be the map such that \(g(s,t)=f(s,t)\) if \(\{s,t\}\neq \{u,v\}\) and \(g(u,v)=g(v,u)^{-1}=f(u,v)r_f(v)^{-1}\), where \(r_f(v)\) is computed by starting from the vertex right after \(u\) in the ordering of \(N_H(v)\). 
Then \(g\) is a (not necessarily tractable) flow in \(H\) such that \(r_g(w)=1\) for \(w\in V\setminus\{u\}\). 
From \(g(v,u)=r_f(v)f(v,u)\) it follows that \(r_g(u)\) differs from \(r_f(u)\) by a factor \(r_f(v)\) when starting the multiplication at \(v\).
By Theorem~\ref{theorem:planar_cycle_stability} we have \(1\equiv r_g(u)\equiv r_f(u)r_f(v)\) and thus \(e_f(u)e_f(v)=1\).
Hence \(G\) is binary leak-proof.

(\(\Rightarrow\))
If \(G=(V,E)\) is not extra-planar, then it has \(K_5^-\) or \(K_{3,3}^-\) as minor by Theorem~\ref{thm:extra_kuratowski}.
It is straightforward to generalize Proposition~\ref{prop:minor_leak} to show that a graph is binary leak-proof if and only if all its minors are too.
It therefore suffices to show that \(K_5^-\) and \(K_{3,3}^-\) have a binary leaking flow.
Simply take the flow \(f\) as defined in Example~\ref{example:k33} which leaks at vertex \(6\) of \(K_{3,3}\) and consider \(K_{3,3}^-\) as the \(K_{3,3}\) with the edge \(\{3,6\}\) removed. Then the flow \(f^-\) of \(K_{3,3}^-\) which equals \(f\) except for \(f^-(3,6)=f^-(6,3)=1\) has a binary leak at \(3\) and \(6\).
Using Example~\ref{example:k5} for \(K_5^-\) can be done analogously. 
\endproof

\section{Leak-proof groups}\label{sec:group}
In this section we prove Theorem~\ref{thm:lpgroup} and give some computational results.

Given a (not necessarily finite) group \(\Gamma\), let \(V(\Gamma)\) be the set of maximal abelian subgroups of \(\Gamma\).
We construct an abelian (additive) group \(\Delta\), which can be defined in various equivalent ways:
\[  \Delta = \text{colim}\bigg( \bigoplus_{(H_1,H_2)\in V(\Gamma)^2 } (H_1\cap H_2)  \doublerightarrow{e_1}{e_2}  \bigoplus_{H\in V(\Gamma)} H\bigg) = \bigg( \bigoplus_{H\in V(\Gamma)} H \bigg) \big/ \im(e_2-e_1) = \coker(e_2-e_1) , \]
where \(e_i\) naturally embeds the summand with index \((H_1,H_2)\) into \(H_i\) and where \(\text{colim}\) is the colimit of this diagram in the category of abelian groups, or equivalently the coequalizer of \(e_1\) and \(e_2\).
%Since the involved groups are abelian, we may explicitly compute \(\Delta\), which is a quotient group of the finite abelian group \(F=\bigoplus_{H\in V(\Gamma)} H\).
For each \(\gamma\in\Gamma\) there is some \(H_\gamma\in V(\Gamma)\) such that \(\gamma\in H_\gamma\) and write \(\overline{\gamma}\) for the image of \(\gamma\) in \(F\) at coordinate \(H_\gamma\).
Now we consider the map \(\varphi_\Gamma:\Gamma\to\Delta\) that naturally sends \(\gamma\) via \(\overline{\gamma}\) to \(\Delta\).
Note that this map does not depend on the choice of \(H_\gamma\):
If we chose \(H_\gamma'\) instead the fact that \(\gamma\in H_\gamma\cap H_\gamma'\) ensures that the image of \(\gamma\) in \(\Delta\) will be the same by definition of \(\Delta\).

\begin{proposition}\label{prop:leak_proof_group}
A group \(\Gamma\) is leak-proof if and only if \(\varphi_\Gamma(\gamma)=0\) implies \(\gamma=1\).
\end{proposition}
\proof
Let \(\gamma\in\Gamma\).
We have \(\varphi_\Gamma(\gamma)=0\) if and only if \(\overline{\gamma}\in\im(e_2-e_1)\), or equivalently there exists an \(x=(x_{HI})_{(H,I)\in V(\Gamma)^2}\in \bigoplus_{(H,I)\in V(\Gamma)^2 } (H\cap I)\) such that \((e_2-e_1)(x)\) equals \(\overline{\gamma}\).
Now consider \(f:V(\Gamma)^2\to\Gamma\) given by \(f(H,I)=x_{HI}\cdot x_{IH}^{-1}\).
Clearly \(f(H,I)=f(I,H)^{-1}\) for all \(H,I\in V(\Gamma)\), while \(\langle f(H,I) \,|\, H\in V(\Gamma)\rangle\subseteq I\) is an abelian group for all \(I\in V(\Gamma)\). 
Hence \(f\) is a tractable flow of \(G=(V(\Gamma),\binom{V(\Gamma)}{2})\). 
Moreover, \((e_f(I))_{I\in V(\Gamma)}=(e_2-e_1)(x)=\overline{\gamma}\), so \(e_f(I)=\gamma\) for \(I=H_\gamma\) and \(e_f(I)=1\) otherwise.
If \(\Gamma\) is leak proof, then \(f\) does not leak and we must have that \(\gamma=1\), proving \((\Rightarrow)\). 
For the converse, note that any flow \(f:V^2\to\Gamma\) in some graph \((V,E)\) gives us a flow for \(G\) by associating to each vertex \(v\in V\) an \(H_v\in V(\Gamma)\) such that \(\langle f(u,v) \,|\, u\in V\rangle \subseteq H_v\) and then applying a contraction.
If \(f\) leaks with value \(\gamma\) at \(v\), then this corresponding flow of \(G\) will also leak with value \(\gamma\) at \(H_v\). 
From the proof above it then follows that \(\varphi_\Gamma(\gamma)=0\), proving \((\Leftarrow)\).
\endproof

Similarly, one can consider binary leak-proof groups. 
With a proof analogous to that of Proposition~\ref{prop:leak_proof_group} one obtains that \(\Gamma\) is binary leak proof if and only if \(\varphi_\Gamma\) is injective.

\proof{Proof of Theorem~\ref{thm:lpgroup}.}
Simply note that for finite \(\Gamma\) the corresponding group \(\Delta\) is finite abelian and can thus be computed explicitly.
In particular, we can decide for each \(\gamma\in\Gamma\) whether \(\varphi(\gamma)=0\). 
The theorem thus follows from Proposition~\ref{prop:leak_proof_group}.
\endproof

From Lemma~\ref{lem:abelian_easy} it follows that abelian groups are leak-proof, but they are hardly the only ones.
By computer search we found the two extraspecial groups of order 32 to be the only smallest leaking groups, one of which we encountered in Example~\ref{example:k33}.
The smallest leaking groups of order greater than 32 occur at order 64.
That there are groups of order 64 that leak was to be expected, because a group leaks when it has a leaking subgroup.
The smallest leaking symmetric group is the \(S_6\) and the smallest leaking alternating group is the \(A_7\).
That for sufficiently large \(n\) the group \(S_n\) leaks is to be expected by Cayley's theorem, but interestingly no strict subgroup of \(S_6\) leaks.
It would be interesting to have a classification of leak-proof groups or to know whether there is some equivalent, better understood property of groups which is equivalent to being leak-proof like planarity is to graphs.

\section*{Acknowledgements}
The author would like to thank H.W. Lenstra for his contributions to Section~\ref{sec:planar_graphs}, for his helpful comments and suggestions and for motivating me to write this article. The author would also like to thank D. Gijswijt for providing references to relevant literature.

\bibliography{citations}{}
\bibliographystyle{plain}
\end{document}